\documentclass[11pt,a4paper]{article}
\usepackage[english]{babel}
\usepackage{amsmath,amssymb,amsfonts,amsthm,latexsym,stmaryrd,mathtools}
\usepackage{multido}
\usepackage{hyperref}
\usepackage{subfig}
\usepackage{graphicx}
\usepackage{xcolor}
\usepackage{mwe}
\usepackage{bm}
\usepackage[round]{natbib}
\usepackage{csquotes}
\usepackage{relsize,etoolbox}
\AtBeginEnvironment{quote}{\smaller}
\usepackage{multirow}
\usepackage{framed}
\usepackage{booktabs}
\usepackage{mathrsfs}
\usepackage{lipsum}
\usepackage{afterpage}
\usepackage[nottoc]{tocbibind}
\usepackage{hyperref}
\def\bkappa1{$\bar{\kappa}_1$}
\def\bkappa2{$\bar{\kappa}_2$}
\def\A{\mathcal{A}_{\varepsilon}}

\def\B{\mathcal{B}_{\varepsilon}}
\def\p{p_0}
\def\pp{p_{\frac{\pi}{2}}}
\def\T{\u{\theta}}

\def\Z{\mathbb{Z}}

\def\R{\mathbb{R}}

\def\Bzero{\mathcal{B}_{\varepsilon,0}}
\def\Bpi{\mathcal{B}_{\varepsilon,\frac{\pi}{2}}}
\def\C{\mathbb{C}}

\def\H{{\mbox H}}

\def\mod{\mbox {\rm mod}}

\def\prob{\stackrel{P}{\longrightarrow}}

\def\e{{\rm e}}

\def\u{\bm}
\def\eps{\varepsilon}

\theoremstyle{plain}
\newtheorem{theorem}{Theorem}[section]

\newtheorem{proposition}[theorem]{Proposition}

\theoremstyle{definition}

\makeindex
\begin{document}
	\newpage
	\vspace*{-5mm}
	
	\centerline{\large BAYESIAN TESTS OF SYMMETRY FOR}
	\vspace{3mm}
	\centerline{\large THE GENERALIZED VON MISES DISTRIBUTION}
	
	\vspace{10mm}
	
	\centerline{Sara Salvador (1) and Riccardo Gatto (2)}
	
	\vspace{6mm}
	
	\centerline{\small Submitted: July 2, 2020}
	
	\centerline{\small Revised: April, 26 2021}
	
	\vspace{8mm}
	
	\centerline{\bf Abstract}
	\noindent 
	Bayesian tests on the symmetry of the generalized von Mises model for planar directions 
	\citep{gatto2007generalized} are introduced. The generalized von Mises distribution 
	is a flexible model that can be axially symmetric or asymmetric, unimodal or bimodal. 
	A characterization of axial symmetry is 
	provided and taken as null hypothesis for one of the proposed Bayesian tests.
	The Bayesian tests are obtained by the technique of probability perturbation.
	The prior probability measure is perturbed so to give  a positive 
	prior probability to the null hypothesis, which would be null otherwise.
	This allows for the derivation of simple computational formulae for the Bayes factors. 
	Numerical results reveal that, whenever 
	the simulation scheme of the samples supports the null hypothesis,
	the null posterior probabilities appear systematically larger than 
	their prior counterpart.
	\vspace{08mm}
	
	\centerline{\bf Key words and phrases}
	\noindent
	{\it
		Axial symmetry;
		Bayes factor;
		circular distribution;
		probability perturbation;
		uni- and bimodality.
	}
	\\\\
	\put(0,0){\line(1,0){462}}
	
	\noindent{\footnotesize
		{\it The authors are grateful to two 
			anonymous Referees and an Associate Editor 
			for several suggestions and corrections.} \\ 
		\\
		{\bf 2010 Mathematics Subject Classification} \\
		\hspace*{1cm} 62H11 Directional data; spatial statistics \\ 
		\hspace*{1cm} 62F15 Bayesian inference \\
		\hspace*{1cm} 62F03 Hypothesis testing \\
		\\
		{\bf Address} \\
		\hspace*{1cm} Institute of Mathematical Statistics and Actuarial Science \\
		\hspace*{1cm} Department of Mathematics and Statistics \\
		\hspace*{1cm} University of Bern \\
		\hspace*{1cm} Alpeneggstrasse 22, 3012 Bern, Switzerland \\
		\hspace*{1cm} (1) {\tt sara.salvador@stat.unibe.ch}, {\tt orcid.org/0000-0001-6492-645X}\\
		\hspace*{1cm} (2) {\tt gatto@stat.unibe.ch}, {\tt orcid.org/0000-0001-8374-6964}
	}
	
	\newpage
	
	\section{Introduction}
	
	In various scientific fields measurements can take the form of directions:
	the direction flight of a bird and the direction of earth's magnetic pole
	are two examples. These directions can be in the plane, namely in two dimensions, 
	as in the first example, or they can be in the space, namely in three dimensions, as 
	in the second example. These measurements are called directional data and they appear
	in various scientific fields:
	in the analysis of protein structure, in machine learning,  
	in forestry, in ornithology, in palaeomagnetism, in oceanography,
	in meteorology, in astronomy, etc. 
	A two-dimensional direction is a point in $\R^2$ without magnitude, e.g. a unit vector.
	It can also be represented as a point on the circumference of the unit circle 
	or as an angle, measured for example in radians and after fixing the null direction
	and the sense of rotation (clockwise or counter-clockwise). Because of this circular representation,
	observations on two-dimensional directional data are distinctively called circular data.
	During the last two or three decades, there has been a  
	raise of interest for statistical 
	methods for directional data.  
	Recent applications can be found e.g. in \cite{ley2018applied}.
	Some monographs on this topic are \cite{mardia2000directional}, \cite{jammalamadaka2001topics},
	\cite{ley2017modern} and also \cite{pewsey2013circular}. For a review article, see e.g.
	\cite{gatto2014directional}.
	
	The popular probability distribution for circular data, 
	or circular distribution, is the circular normal or von Mises 
	distribution, whose density is given in (\ref{e6}) below.
	This distribution is circularly symmetric around its unique mode. 
	Until a couple decades ago, very few asymmetric circular distributions 
	were available, two of these can be found in 
	Sections 15.6 and 15.7 of 
	\cite{batschelet1981circular}.
	In recent years, various asymmetric and multimodal circular distributions have been
	introduced, for example: \cite{umbach2009building}, 
	\cite{kato2015tractable}, \cite{abe2013extending}, \cite{gatto2003inference} 
	and the generalized von Mises (GvM) of \cite{gatto2007generalized}.
	This article proposes three Bayesian tests for the 
	GvM distribution. This distribution has density given  
	\begin{align}  							                      \label{e4}
		f(\theta \mid \mu_1, \mu_2, \kappa_1, \kappa_2) & =
		\frac{1}{2 \pi G_0(\delta,\kappa_1,\kappa_2)}
		\exp \{\kappa_1 \cos(\theta -\mu_1) + \kappa_2 \cos 2(\theta-\mu_2)\},
	\end{align}
	$\forall \theta \in [0, 2 \pi)$, for given $\mu_1 \in [0, 2 \pi)$, $\mu_2 \in [0,\pi)$,
	$\delta = (\mu_1 - \mu_2) \mod \, \pi$, $\kappa_1,\kappa_2>0$,
	and where the normalizing constant is given by
	\begin{align}              				  			        \label{e5}
		G_0(\delta,\kappa_1,\kappa_2) & =
		\frac{1}{2 \pi} \int_0^{2 \pi} \exp\{\kappa_1 \cos \theta + \kappa_2 \cos 2(\theta+\delta)\}
		d \theta.
	\end{align}
	We denote any circular random variable with this distribution by 
	GvM($\mu_1,\mu_2, \kappa_1,\kappa_2$).
	The well-known von Mises (vM) density is obtained by setting $\kappa_2 = 0$ in
	(\ref{e4}), giving
	\begin{align}           			                                     \label{e6}
		f(\theta \mid \mu,\kappa) & =
		\frac{1}{2 \pi I_0(\kappa)} \exp \{\kappa \cos(\theta -\mu)\},
	\end{align}
	$\forall \theta \in [0, 2 \pi)$, for given $\mu \in [0, 2 \pi)$, $\kappa>0$ and
	where $I_{\nu}(z)= (2 \pi)^{-1} \int_0^{2 \pi} \cos \nu \theta \, \exp\{
	z \cos \theta\} d \theta $, $z \in \C$, 
	is the modified Bessel function $I$ of order $\nu$, with $\Re \nu > -1/2$, 
	cf. 9.6.18 at p. 376 of \cite{abramowitz1972handbook}.  
	We denote any circular random variable $\theta$ with this distribution by
	vM($\mu,\kappa$).
	
	Besides its greater flexibility in terms of asymmetry and bimodality,
	the GvM distribution possesses the following properties that other
	asymmetric or multimodal circular distributions do not have. 
	\begin{enumerate} 
		\item  After a re-parametrization, the GvM distribution belongs to the canonical exponential class. In this form,
		it admits a minimal sufficient and complete statistic; cf. Section 2.1 of 
		\cite{gatto2007generalized}. 
		\item
		The maximum likelihood estimator and the trigonometric method of moments estimator of the parameters
		are the same; cf. Section 2.1 of \cite{gatto2008some}.  
		In this context, we should note that
		the computation of the maximum likelihood estimator 
		is simpler with the GvM distribution than with the mixture of two vM distributions, 
		as explained some lines below. 
		\item It is shown in Section 2.2 of \cite{gatto2007generalized} that
		for fixed trigonometric moments of orders one and two, the GvM distribution is the one with 
		largest entropy. The entropy gives a principle for selecting a distribution 
		on the basis of partial knowledge: one should always choose distributions having maximal entropy, 
		within distributions satisfying the partial knowledge.
		In Bayesian statistics, 
		whenever a prior distribution has to be selected and
		information on the first two trigonometric moments is available, then the GvM is the optimal prior.
		For other theoretic properties of the GvM, see \cite{gatto2009information}. 
	\end{enumerate}
	The mixture of two vM distributions is perhaps a more popular bimodal or asymmetric
	model then the GvM. However, the mixture does not share the 
	given properties 1-3 of the GvM. The mixture is not necessarily more practical. 
	While the likelihood of the GvM distribution is bounded, the
	likelihood of the mixture of the vM($\mu_1,\kappa_1$) and the vM($\mu_2,\kappa_2$)
	distributions is unbounded. As $\kappa_1 \rightarrow \infty$, the likelihood when 
	$\mu_1$ is equal to any one of the sample values tends to infinity. This follows from
	$I_0(\kappa_1) \sim ( 2 \pi \kappa_1)^{-1/2} \e^{\kappa_1}$, 
	as $\kappa_1 \rightarrow \infty$; cf. \cite{abramowitz1972handbook}, 9.7.1 at p. 377. For
	alternative estimators to the maximum likelihood for vM mixtures, refer to
	\cite{spurr1991comparison}.
	
	Some recent applications of the GvM distributions are:
	\cite{Zhang18}, in meteorology, \cite{Lin19}, in oceanography,
	\cite{astfalck2018expert}, in offshore engineering, 
	\cite{6784494}, in signal processing, and 
	\cite{gatto2021spectrum} in time series analysis.
	
	The symmetry of a circular distribution is a fundamental question and, as 
	previously mentioned, this topic has been studied in recent years.
	In the context of testing symmetry, one can mention: 
	\cite{pew02}, who proposes a test of symmetry around an unknown axis 
	based on the second sine sample moment, and 
	\cite{pew04}, who considers the case where the symmetry is
	around the median axis. Both tests are frequentist and no Bayesian 
	test of symmetry appears available in the literature.
	In fact, Bayesian analysis for circular data has remained underdeveloped,
	partly because of the lack of nice conjugate classes of distributions. 
	Moreover, Bayesian analysis has focused on the vM model, which is symmetric.
	We refer to p. 278-279 of \cite{jammalamadaka2001topics} 
	for a review on Bayesian analysis for circular data. 
	
	In this context,
	this article proposes Bayesian tests of symmetry for the GvM model  
	(\ref{e4}). 
	The first test proposed concerns the parameter $\delta$. The null hypothesis
	is $\delta = 0$, that is, no shift between cosines of frequency 
	one and two. In this case, the distribution is symmetric around the axis passing
	through $\mu_1$. It is bimodal with one mode at $\mu_1$ and the other one
	at $\mu_1 + \pi$, whenever $\kappa_1 < 4 \kappa_2$.
	If $\kappa_1 \geq 4 \kappa_2$, then it is unimodal with mode at $\mu_1$. We refer
	to Table 1 of \cite{gatto2007generalized}.   
	The second test is on the precise characterization of axial symmetry, i.e.  
	on $\delta = 0$ or $\delta = \pi/2$. 
	So far $\kappa_2 >0$ is considered and the third test is for $\kappa_2 = 0$, so that
	the distribution is no longer GvM but vM, which is 
	is axially symmetric. 
	The Bayesian tests rely on the method of probability perturbation, where the probability 
	distribution of the null hypothesis is slightly perturbed,
	in order to give a positive prior probability to the null hypothesis, which
	would be null otherwise. 
	It would be interesting to consider 
	the above null hypotheses under the frequentist perspective, perhaps with the 
	likelihood ratio approach. 
	This topic is not studied in this article, in order to limit its length. 
	
	The remaining part of this article is organized as follows. 
	Section \ref{s2} gives the derivation of these Bayesian tests and 
	their Bayes factors. Section \ref{s2new} presents the approach used
	for these tests:
	Section \ref{s21} considers the test
	of no shift between cosines, Section \ref{s22} considers
	the test of symmetry and Section \ref{s23} considers 
	the test of vM axial symmetry.
	Numerical results are presented in Section \ref{s3}:
	Section \ref{s310} presents a Monte Carlo study of the the tests of 
	Section \ref{s2new} whereas
	Section \ref{s40} presents the application to some real data.
	Final remarks are given in Section \ref{s4}. 
	
	\section{Bayesian tests and perturbation method for the GvM model}	\label{s2}
	
	The proposed tests rely on Bayes factors. The Bayes factor
	$B_{01}$ indicates the evidence of the null hypothesis with respect to (w.r.t.) 
	the general 
	alternative.
	Let us denote by $\u{\theta} = ( \theta_1,\ldots,\theta_n )$ the sample. Then
	\begin{equation}
		B_{01}=\frac{P[\u\theta|\H_0]}{P[\u\theta|\H_1]}=\frac{P[\H_0|\u\theta]}
		{P[\H_1|\u\theta]}\cdot \frac{P[\H_1]}{P[\H_0]}=\frac{R_1}{R_0},
		\label{BF}
	\end{equation}
	where
	\[
	R_0=\frac{P[\H_0]}{P[\H_1]}=\frac{P[\H_0]}{1-P[\H_0]}~~\text{and}~~ R_1=\frac{P[\H_0|
		\u\theta]}{P[\H_1|\u\theta]}= \frac{P[\H_0|\u\theta]}{1-P[\H_0|\u\theta]},
	\]
	are the prior and the posterior odds, respectively. 
	The case
	$B_{01}>1$ indicates evidence for $\H_0$. Interpretations of
	the values of the Bayes factor can be found in 
	\cite{jeffreys1961theory} and \cite{kass1995bayes}. 
	Our synthesis of these interpretations is given in Table \ref{tab}, which 
	provides a qualitative scale for the Bayes factor.    
	\begin{table}
		{\small 
			\centering
			\begin{tabular}{|c|c|}
				\hline
				$B_{01}$ & evidence for $\H_{0}$\\
				\hline
				$<1$ & negative \\
				1 to 1.5 &  not worth more than a bare mention\\
				1.5 to 5 & positive \\
				5 to 10 & substantial\\
				10 to 20 & strong\\
				$>20$ & decisive\\
				\hline
			\end{tabular}
			\caption{\small Guidelines for the interpretation of Bayes factors.}						\label{tab}
		}
	\end{table}
	The null hypotheses of this article are simple, in the sense that they concern
	only points of the parametric space. The fact that these points have probability null
	does not allow for the computation of Bayes factors.
	Therefore we use an approach with probability perturbation explained
	in the next section.
	
	\subsection{Bayesian tests of simple hypotheses}		\label{s2new} 
	
	The practical relevance of a simple null hypothesis, i.e. of the type 
	$\H_0:\xi=\xi_0$, has been widely 
	debated in the statistical literature.  
	According to
	\citeauthor{bergerdelampady1987testing}:  
	``it is rare, and perhaps impossible, to have a null hypothesis that 
	can be exactly modelled as $\theta=\theta_0$''.
	They illustrate their claim by the following example.
	``More common precise hypotheses such as $\H_0$:\emph{Vitamin C has no effect on the common cold} are clearly not meant to be though of as exact point nulls; surely vitamin C has some effects, although perhaps a very miniscule effect.'' 
	A similar example involving forensic science can be found in \cite{lindley1977problem}. 
	When the parameter $\xi$ is of continuous nature, 
	it is usually more realistic to consider  
	null hypotheses of the type $\H_{0,\eps}:|\xi - \xi_0|\leq \varepsilon/2$, 
	for some small $\varepsilon>0$. This solves also the problem of the vanishing 
	prior probability of $\H_0$, namely $P[\xi=\xi_0]=0$. This problem is sometimes 
	addressed by giving a positive probability to $\{\xi=\xi_0\}$. However,
	\cite{berger1987testing} explain that the two approaches should be related.
	``It is convenient to specify a prior distribution for the testing problem as follows: let $0<\pi_0<1$ denote the prior probability of $\H_0:\theta=\theta_0$ ... One might question the assignment of a positive probability to $\H_0$, because it is rarely 
	the case that it is thought possible for $\theta=\theta_0$ to hold exactly ... 
	$\H_0$ is to be understood as simply an approximation to the realistic hypothesis 
	$\H_0:|\theta-\theta_0|\leq b$ and $\pi_0$ is to be interpreted as the prior probability that would be assigned to $\{\theta:|\theta-\theta_0|\leq b\}$.''
	Accordingly, we assign to the original simple hypothesis $\H_0: \xi = \xi_0$ 
	the prior probability $p_0>0$ of 
	$\H_{0,\eps}:\xi\in[\xi_0-\varepsilon/2, \xi_0+\varepsilon/2]$, for some $\eps>0$.
	Thus, we replace the prior probability measure $P$ 
	by its perturbation, obtained by the assignment of the probability $p_0>0$
	to $\{\xi_0\}$.
	We denote by $P_0$ the probability measure $P$
	with the $p_0$-perturbation.
To summarize: the point null hypotheses is made relevant with 
	$p_0=P_0[\xi=\xi_0]=
	P\left[\delta\in\left[\xi_0-\varepsilon/2, \xi_0+\varepsilon/2\right]\right]>0$.
	
	The length  
	$\eps$ of the neighbourhood of $\xi_0$, which determines the prior probability 
	$p_0$ of $\H_0$ under the perturbed model,
	should not be too small.
	A significant value of $p_0$ for  
	the null hypothesis is in fact coherent with the frequentist 
	approach of hypotheses tests, where
	computations of rejection regions or P-values are carried over 
	under the null hypothesis. 
	\cite{berger1985statistical}, p. 149, states that
	$\varepsilon$ has to be chosen such that any $\xi$ in 
	$(\xi_0-\eps/2,\xi_0+\eps/2)$ becomes ``indistinguishable'' from $\xi_0$, 
	while \citeauthor{berger1987testing} state that $\varepsilon$ 
	has to be ``small enough'' so that $\H_{0,\eps}$ can be 
	``accurately approximated'' by $\H_0$. A related reference is 
	\cite{bergerdelampady1987testing}, who studied this problem 
	with a Gaussian model, and \cite{berger1985statistical}, p. 149, who obtains an upper 
	bound for the radius $\varepsilon/2$ under a simple Gaussian model. 
	Two other references on the practical relevance 
	of simple null hypotheses are \cite{jeffreys1961theory} and \cite{Zellner1984}. 
	
	We end this section with some comments regarding the choice of the 
	prior distribution of $\xi$. This is a generally 
	unsolved problem of Bayesian statistics and widely discussed in the 
	literature, see e.g. \cite{jeffreys1961theory} and \cite{kass1996selection}. 
	According to \cite{bergerdelampady1987testing}, there is
	``no choice of the prior that can claim to be objective''. 
	In this article we follow the directives given in \cite{bergerdelampady1987testing} and \cite{berger1987testing}, where various details on the choice of the prior are discussed and some classes of priors are analysed. According to \cite{bergerdelampady1987testing}, 
	in absence of prior information, the prior should be 
	symmetric about $\xi_0$ and non-increasing w.r.t. $|\xi-\xi_0|$. Otherwise, 
	one could find some ``favoured'' alternative values of $\xi$; 
	cf. \cite{berger1987testing}. Our choices 
	of priors are presented in Section \ref{s3}: for each test of the study 
	we compute Bayes factors under priors obtained by
	varying the concentration around the generic value $\xi_0$. 

	\subsection{Test of no shift between cosines of GvM}	\label{s21}
	
	Consider the Bayesian test on the GvM model (\ref{e4})
	of the null the hypothesis
	\[
	\text{H}_{0}: \delta=0,
	\] 
	where $\delta=(\mu_1-\mu_2)~\text{mod}~\pi$
	and where the values of $\mu_1,\kappa_1,\kappa_2$ are assumed known and equal to
	$\mu_1^0,\kappa_1^0,\kappa_2^0$, respectively. 
	Under the original probability measure $P$, the random parameter $\delta$ 
	has an absolutely continuous prior distribution and so
	$P[\delta=0]=0$. 
	According to Section \ref{s2new}
	we define the perturbation of the probability measure $P$, denoted $P_0$, 
	for which
	$p_0=P_0[\delta=0]>0$.
	This perturbation is the assignment to $\{\delta=0\}$ 
	of the probability mass that initially lies close to that $P$-null set. 
	Let $\eps>0$ and consider the set
	\begin{align}								\label{e13} 
		\mathcal{A}_{\varepsilon} & = \left\{\delta\in[0,\pi) \Big|  \delta\in \left[0,\frac{\varepsilon}{2}\right]
		\cup\left[\pi-\frac{\varepsilon}{2},\pi\right)\right\}.
	\end{align}
	The complement is
	\begin{align*}
		\mathcal{A}^{{\sf c}}_\eps & = \left\{\delta\in[0,\pi) \Big| \delta\in \left(\frac{\varepsilon}{2},\pi-\frac{\varepsilon}{2}\right)\right\}.
	\end{align*}
       Note that (\ref{e13}) refers to a neighbourhood of the origin of the circle of circumference
        $\pi$.
	We thus assign to $p_0$ the value  
	\begin{align}                                                                                           \label{e10}
		p_0 & = P\left[ \mathcal{A}_{\varepsilon} \right],
	\end{align}
	for some suitably small $\varepsilon>0$.
	The prior distribution function (d.f.) 
	under the perturbed probability measure $P_0$ at any $\delta \in [0,\pi)$ 
	is given by
	\begin{align}												\label{e23} 
		p_0 \Delta ( \delta' ) + (1 - p_0 ) G ( \delta'). 
	\end{align}
	where $G$ denotes the prior d.f. of $\delta$ and where $\Delta$ is the Dirac d.f., which assigns mass one to the 
	origin. Denote by $g$ the density of $G$.
	If $0\notin (\delta',\delta'+d\delta')$, for some $\delta' \in (0, \pi)$, where
	the relations $\in$ and $\notin$ are meant circularly over the circle of circumference $\pi$, then
	(\ref{e23}) implies 
	\begin{equation}						\label{e34}
		\begin{split}
			P_0[\delta\in(\delta',\delta'+d\delta')] & =  (1 - p_0) g (\delta') d \delta' = (1 - p_0) P[\delta\in(\delta',\delta'+d\delta')].
		\end{split}
	\end{equation}
	
	Let $\theta_1,\ldots,\theta_n$ be independent circular random variables that follow the GvM distribution 
	(\ref{e4}). 
	For simplicity, we denote the joint density of $\u{\theta}=(\theta_1,\ldots,\theta_n)$, with the
	fixed values $\delta'$, $\mu_1^0$, $\kappa_1^0$ and $\kappa_2^0$, as
	\begin{align}                                                                                   \label{dd}
		f(\u{\theta}|\delta') & = 
		\left\{2\pi G_0(\delta',\kappa_1^0,\kappa_2^0)\right\}^{-n}
		\exp\left\{\kappa_1^0 \sum_{i=1}^{n}\cos(\theta_i-\mu_1^0)
		+\kappa_2^0 \sum_{i=1}^{n} \cos 2 ( \theta_i - \mu_1^0 +\delta' ) \right\}.
	\end{align}
	When considered as a function of $\delta'$, (\ref{dd}) becomes the likelihood of $\delta$.
	Then, by (\ref{e34}) the marginal density of $\u{\theta}=(\theta_1,\ldots,\theta_n)$ under the perturbed probability 
	is given by
	\begin{align}								\label{e278}	
		m(\u{\theta})&=\int_{[0,\pi)} f(\u{\theta}|\delta')
		P_0[\delta\in (\delta',\delta'+d\delta')] \nonumber \\
		&=\int_{\mathcal{A}_{\varepsilon}} f(\u{\theta}|\delta')P_0[\delta\in(\delta',\delta'+d\delta')]+(1-p_0)\int_{\mathcal{A}_{\varepsilon}^{\mathsf{c}}}f(\u{\theta}|\delta')g(\delta')d\delta' 
		\nonumber \\
		& = p_0f(\u{\theta}|0)
		+(1-p_0)\int_{\mathcal{A}_{\varepsilon}}f(\u{\theta}|\delta')
		g(\delta')d\delta'
		+(1-p_0)\int_{\mathcal{A}_{\varepsilon}^{\mathsf{c}}}f(\u{\theta}|\delta')
		g(\delta')d\delta' \nonumber \\
		&\sim 2 p_0f(\u{\theta}|0)+(1-p_0)\int_{\mathcal{A}_{\varepsilon}^{\mathsf{c}}}f(\u{\theta}|\delta')g(\delta')d\delta',~\text{as $\varepsilon\rightarrow0$}.
	\end{align}
	The above asymptotic equivalence is due to 
	\begin{align*}
		(1-p_0)\int_{\mathcal{A}_{\varepsilon}}f(\u{\theta}|\delta') g(\delta')d\delta' 
		& = (1 - p_0) p_0 \int_{\mathcal{A}_{\varepsilon}}f(\u{\theta}|\delta') \frac{g(\delta')}{p_0} d\delta' 
		\sim p_0 f( \theta | 0 ), \; \text{as $\varepsilon\rightarrow0$}.
	\end{align*}
	
	The posterior perturbed probability, 
	namely the conditional perturbed probability 
	of $\{\delta = 0 \}$ given $\u{\theta}$, 
	can be approximated as follows,
	\begin{equation}
		\begin{split}
			P_0[\delta=0|\u{\theta}]
			&\sim\frac{ p_0f(\u{\theta}|0)}{2 p_0 f(\u{\theta}|0)+(1-p_0)\int_{\A^\mathsf{c}}f(\u{\theta}|\delta')
				g(\delta')d\delta'}  \\
			&=\frac{1}{2}\Bigg(1+\frac{1-p_0}{p_0}\frac{\int_{\A^\mathsf{c}}f(\u{\theta}|\delta')g(\delta')d\delta'}{2 f(\u{\theta}|0)}\Bigg)^{-1}, \;
			\text{ as } \varepsilon \rightarrow 0.
		\end{split}
		\label{post}
	\end{equation}
	
	In order to compute the Bayes factor for this test, we define the prior odds  
	$R_0 = p_0/(1-p_0)$ and
	the posterior odds  
	$R_1 = P_0[\delta=0|\u{\theta}]/(1-P_0[\delta=0|\u{\theta}])$.
	The Bayes factor is the posterior over the prior odds, namely 
	$
	B_{01}=R_1/R_0$.
	Clearly $p_0 \le P_0[\delta=0|\u{\theta}]$ iff $B_{01} \ge 1$ and, the 
	larger $P_0[\delta=0|\u{\theta}] - p_0$ becomes, the larger $B_{01}$ becomes: a large Bayes factor tells that the 
	data support the null hypothesis. From the approximation
	\begin{equation*}
		\begin{split}
			R_1&\sim\left[1+\frac{1-p_0}{p_0}\frac{\int_{\A^{\mathsf{c}}}
					f(\u{\theta}|\delta')g(\delta')d\delta'}{f\left(\T|0\right)}\right]^{-1}
		\end{split}
	\end{equation*}
	and from some simple algebraic manipulation, we obtain the computable approximation to the Bayes factor $B_{01} = R_1/R_0$ given by
	\begin{equation}			
		B_{01}
		\sim\frac{ f(\u{\theta}|0)}{\int_{\A^{\mathsf{c}}}f(\u{\theta}|\delta')g(\delta')d\delta'},~\text{as $\varepsilon\rightarrow0$}.
		\label{BF_delta}
	\end{equation}
	The representation of the Bayes factor \eqref{BF_delta} is asymptotically correct and we remind that, in the context where we approximate the null hypothesis with a neighbourhood by the point null hypothesis, the reasoning is always of asymptotic nature.
	A reference for this perturbation technique is \cite{berger1985statistical}, p. 148-150. 
	
	\begin{color}{black}
		Regarding the large sample asymptotics of the proposed test,
		it is know that, for a sample of $n$ independent random variables with common
		distribution with true parameter $\xi_0$,
		the posterior distribution converges to the distribution with total mass over $\xi_0$,
		as $n \to \infty$. This means that the posterior mode is a consistent estimator. 
		We deduce that, under $\H_0$,
		\begin{align*}
			P_0 [ \delta = 0 | \u{\theta} ] & = P [ {\cal A}_\eps | \u{\theta} ] \prob 1, \;
			\text{ as } n \to \infty.
		\end{align*}
		Consequently, 
		$R_1 = P_0[\delta=0|\u{\theta}]/(1-P_0[\delta=0|\u{\theta}]) \prob \infty$ 
		and $B_{01} = R_1/R_0 \prob \infty$, as
		$n \to \infty$. The Bayesian test of $\H_0: \delta = 0$ is consistent in this sense. 
	\end{color}
	
	We now give some computational remarks that are also valid for the tests 
	of Sections \ref{s22} and \ref{s23}. The integral appearing in the denominator 
	of (\ref{BF_delta}) can be easily evaluated by Monte Carlo
	integration. For a given large integer $s$, 
	we generate $\delta^{{(i)}}$, for $i = 1,\ldots,s$, from the density $g$ 
	and then we compute the approximation
	\begin{align}							\label{e123}
		\int_{\A^{\mathsf{c}}}f(\u{\theta}|\delta')g(\delta')d\delta' 
		& = \int_{\A^{\mathsf{c}}}f(\u{\theta}|\delta') P[ \delta \in (\delta' , \delta' + d \delta') ] \simeq \frac{1}{s}\sum_{i=1}^{s}
		f(\u{\theta}|\delta^{(i)})\text{I}\{\delta^{(i)}\in\A^{\mathsf{c}}\},
	\end{align}
	where $\text{I}\{A\}$ denotes the indicator of statement or event $A$.
	For the computation normalizing constant of the GvM distribution given in (\ref{e5})
	one can use the Fourier series 
	\begin{equation}				\label{e120}        		
		G_0(\delta,\kappa_1,\kappa_2)=
		I_0(\kappa_1) I_0(\kappa_2)
		+ 2 \sum_{j=1}^\infty I_{2j}(\kappa_1)
		I_j(\kappa_2) \cos 2 j \delta,
	\end{equation}
	where $\delta \in [0, \pi)$ and $\kappa_1,\kappa_2 > 0$,
	see e.g. \cite{gatto2007generalized}. 
	\subsection{Test of axial symmetry of GvM}			\label{s22}
	
	In this section we consider the Bayesian test of axial symmetry for the GvM model (\ref{e4}).
	A circular density $g$ is symmetric around the angle 
	$\alpha/2$, for some $\alpha \in [0 , 2 \pi)$, if $g(\theta) = g(\alpha-\theta)$, $\forall \theta 
	\in [0, 2\pi)$. In this case we have also $g( \theta ) = g( (\alpha+2\pi)-\theta)$, so that
	symmetry around $\alpha/2+\pi$ holds as well: the symmetry is indeed an axial one.
	\begin{proposition}[Characterization of axial symmetry for the GvM distribution]		\label{p1}
		The GvM distribution (\ref{e4}) is axial symmetric iff
		$$\delta=0 \; \text{ or } \; \delta=\frac{\pi}{2}.$$
		In both cases, the axis of symmetry has angle $\mu_1$.
	\end{proposition}
	
	The proof of Proposition \ref{p1} is given in Appendix \ref{a1}.
	
	Note that $\delta$ is defined modulo $\pi$ and that
	for $\kappa_2=0$ or $\kappa_1=0$
	the GvM reduces respectively to the vM or to the axial vM, defined later as 
	$\text{vM}_2$ and given in (\ref{e321}).
	These two distributions are clearly symmetric,
	but Proposition \ref{p1} gives the characterization of symmetry 
	in terms of $\delta$ since we define the GvM distribution in (\ref{e4}) 
	with concentration parameters $\kappa_1,\kappa_2>0$.
	
\begin{figure} [h!]
		\centering
		\subfloat[GvM($\pi,0,0.1,5.5$) density ($\delta = 0$).]{\includegraphics[scale=0.3]{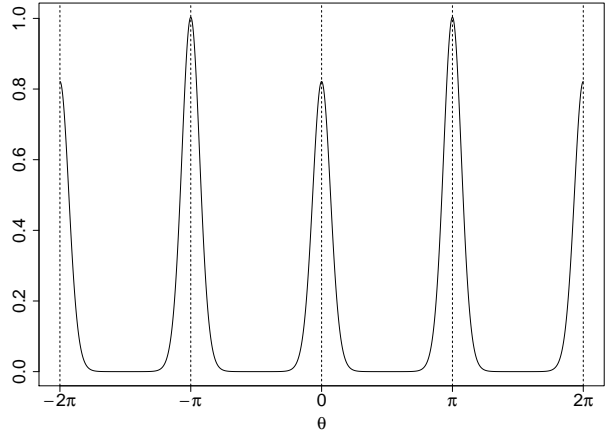}{\label{symm1}}}\quad
		\subfloat[GvM($\pi/2,0,5.5,0.1$) density ($\delta = \pi/2$).]{\includegraphics[scale=0.3]{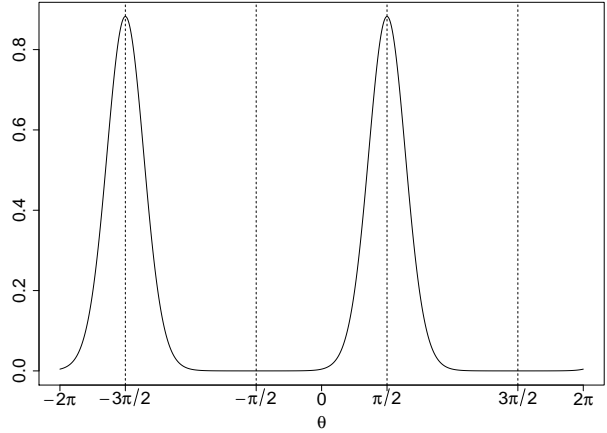}\label{symm2}}\\
		\caption{Two axial symmetric GvM densities 
			over the interval $(-2\pi,2\pi)$ and with their axis of symmetry at angle $\mu_1$ 
			shown by vertical dashed lines.}
		\label{ax_sim}
	\end{figure}

	As mentioned at the beginning of the section, symmetry of a circular distribution 
	around an angle is the symmetry around an axis. For the GvM density,
	this is made explicit in \eqref{sim}, where adding $2 \pi$ to $\alpha$ would not have any influence.
	Figure \ref{ax_sim} provides two numerical illustrations 
	of the axial symmetry of the GvM distribution.
	The graph in Figure \ref{symm1} shows the density of the GvM($\pi,\pi,0.1,5.5$) distribution: $\delta = 0$ and the 
	axis of symmetry is at angle $\mu_1 = \pi$.
	The graph in Figure \ref{symm2} shows the density of the GvM($\pi/2,0,5.5,0.1$) distribution: $\delta = \pi/2$ and the 
	axis of symmetry is at angle $\mu_1 = \pi/2$. 
	Thus, Proposition \ref{p1} allows us to write the null hypothesis of axial symmetry as
	\[
	\text{H}_0: \delta=0  \text{ or } \delta=\frac{\pi}{2},
	\]
	where the values of $\mu_1,\kappa_1,\kappa_2$ are assumed known and equal to
	$\mu_1^0,\kappa_1^0,\kappa_2^0$, respectively. 
	The Bayesian test is obtained by perturbation of the probability measure $P$, which is denoted $P_0$. 
	The probabilities
	\begin{equation*}
		p_0=P_0[\delta=0] > 0 \; \text{ and } \;
		p_{\frac{\pi}{2}}=P_0 \left[ \delta=\frac{\pi}{2} \right] >0
	\end{equation*}
	are the probabilities masses of $\{\delta=0\}$ and $\{\delta=\frac{\pi}{2}\}$ of the perturbed measure, respectively. 
	They are
	obtained from  
	\begin{equation*}
		p_0=P\left[\delta\in\left[0,\frac{\varepsilon}{2}\right]\cup\left[\pi-\frac{\varepsilon}{2},\pi\right)\right] \; \text{ and }
		p_{\frac{\pi}{2}}=P\left[\delta\in\left[\frac{\pi}{2}-\frac{\varepsilon}{2},\frac{\pi}{2}+\frac{\varepsilon}{2}\right]\right],
	\end{equation*}
	for suitably small $\eps >0$.
	As is Section \ref{s21}, the prior d.f. of $\delta$ under 
	the perturbed probability $P_0$ at any $\delta' \in [0,\pi)$ is given by
	\begin{equation}
		p_0\Delta(\delta')+p_{\frac{\pi}{2}}\Delta\left(\delta'-\frac{\pi}{2}\right)
		+\left\{1-\left(p_0+p_{\frac{\pi}{2}}\right)\right\} G(\delta'),
		\label{dfsymm}
	\end{equation}
	where $G$ is the prior d.f. of $\delta$ under $P$.
	It follows from \eqref{dfsymm} that for $0,\pi/2 \notin (\delta',\delta'+d\delta')$, for some $\delta'\in(0,\pi)\setminus\{ \pi / 2 \}$,
	\begin{equation*}
		P_0[\delta\in(\delta', \delta'+d\delta')]=\left[1-\left(p_0+p_{\frac{\pi}{2}}\right)\right]g(\delta')d\delta'=\left[1-\left(p_0+p_{\frac{\pi}{2}}\right)\right]P[\delta\in(\delta',\delta'+d\delta')],
	\end{equation*}
	where $g$ is the density of $G$.
	
		Let 
	\[
	\Bzero=\left\{\delta\in[0,\pi)\Big|\delta\in\left[0,\frac{\varepsilon}{2}\right]\cup\left[\pi-\frac{\varepsilon}{2},\pi\right)\right\},~\text{ and }
	\Bpi=\left\{\delta\in[0,\pi)\Big| \delta\in\left[\frac{\pi}{2}-\frac{\varepsilon}{2},\frac{\pi}{2}+\frac{\varepsilon}{2}\right]\right\}.
	\]
	Define 
	\begin{equation*}
		\B=\left\{\delta\in[0,\pi) \Big| \delta\in\left[0,\frac{\varepsilon}{2}\right]\cup\left[\pi-\frac{\varepsilon}{2},\pi\right)\vee\delta\in\left[\frac{\pi}{2}-\frac{\varepsilon}{2},\frac{\pi}{2}+\frac{\varepsilon}{2}\right]\right\}=\Bzero\cup\Bpi.
	\end{equation*}
	Its complement is given by
	\begin{equation*}
		\B^{\mathsf{c}}=\left\{\delta\in[0,\pi) \Big|  \delta\in\left(\frac{\varepsilon}{2},\frac{\pi}{2}-\frac{\varepsilon}{2}\right)\cup\left(\frac{\pi}{2}+\frac{\varepsilon}{2},\pi-\frac{\varepsilon}{2}\right)\right\}.
	\end{equation*}
	The marginal density of $\u{\theta}=(\theta_1,\ldots,\theta_n)$ w.r.t. the perturbed probability $P_0$ is given by
		\begin{equation*}
		\begin{split}
			m(\T)&=\int_{[0,\pi)} f(\T|\delta')P_0[\delta\in(\delta',\delta'+d\delta')] 
			\nonumber \\
			& =       \int_{\B} f(\T|\delta')P_0[\delta\in(\delta',\delta'+d\delta')]
			+\left[1-\left(p_0+p_{\frac{\pi}{2}}\right)\right]
			\int_{\B^{\mathsf{c}}} f(\T|\delta')g(\delta')d\delta' \nonumber \\
			&=p_0f(\T|0)+p_{\frac{\pi}{2}}f\left(\T|\frac{\pi}{2}\right)+\left[1-\left(p_0+p_{\frac{\pi}{2}}\right)\right]\int_{\B} f(\T|\delta')g(\delta')d\delta'\\
			&\quad\qquad+\left[1-\left(p_0+p_{\frac{\pi}{2}}\right)\right]\int_{\B^{\mathsf{c}}} f(\T|\delta')g(\delta')d\delta'\\
			& \sim    2p_0f(\T|0)+2p_{\frac{\pi}{2}}f\left(\T \Big| \frac{\pi}{2}\right)+\left[1-\left(p_0+p_{\frac{\pi}{2}}\right)\right]
			\int_{\B^{\mathsf{c}}} f(\T|\delta')g(\delta')d\delta',~\text{as $\varepsilon\rightarrow0$}.
		\end{split}
	\end{equation*}
	In the asymptotic equivalence, as in Section \ref{s21}, we notice that
	\begin{equation*}
		\begin{split}
			&\left[1-\left(p_0+p_{\frac{\pi}{2}}\right)\right]\int_{\B} f(\T|\delta')g(\delta')d\delta'\\	&=\left[1-\left(p_0+p_{\frac{\pi}{2}}\right)\right]\left[\int_{\Bzero} f(\T|\delta')g(\delta')d\delta'
			+\int_{\Bpi} f(\T|\delta')g(\delta')d\delta'\right]\\
			&=\left[1-\left(p_0+p_{\frac{\pi}{2}}\right)\right]\left[p_0\int_{\Bzero} f(\T|\delta')\frac{g(\delta')}{p_0}d\delta'+p_{\frac{\pi}{2}}\int_{\Bpi} f(\T|\delta')\frac{g(\delta')}{p_{\frac{\pi}{2}}}d\delta'\right]\\
			&\sim p_0f(\T|0)+p_{\frac{\pi}{2}}f\left(\T\Big|\frac{\pi}{2}\right),~\text{ as }\varepsilon\rightarrow 0.
		\end{split}
	\end{equation*}
	The posterior probability of 
	$\{\delta=0\vee\delta=\pi/2\}$ under the perturbed probability measure is given by
	\begin{equation*}
		\begin{split}
			P_0\left[\delta=0\vee\delta=\frac{\pi}{2} \Big| \T \right]&
			\sim \frac{p_0f(\T\lvert 0)+p_{\frac{\pi}{2}}f(\T|\frac{\pi}{2})}{2p_0f(\T|0)+2p_{\frac{\pi}{2}}f(\T|\frac{\pi}{2})+
				[1-(p_0+p_{\frac{\pi}{2}})]\mathcal{I}_1(\theta)} \\
			&=\frac{1}{2}\left[1+\frac{\left[1-(p_0+p_{\frac{\pi}{2}})\right]\mathcal{I}_1(\u{\theta})}{2p_0f(\T|0)+2p_{\frac{\pi}{2}}f(\T|\frac{\pi}{2})}\right]^{-1},
			\; \text{ as } \eps \to 0,
		\end{split}
	\end{equation*}
	where 
	\begin{equation*}
		\mathcal{I}_1 (\u{\theta}) = \int_{\B^{\mathsf{c}}}
		f(\T|\delta')g(\delta')d\delta'.
		\label{approx}
	\end{equation*}
	With this we obtain the following approximation to the posterior odds,
	\begin{equation*}
		\begin{split}
			R_1=\frac{P_0\left[\delta=0\vee\delta=\frac{\pi}{2}|\T\right]}{1-P_0\left[\delta=0\vee\delta=\frac{\pi}{2}|\T\right]}=\left[\frac{1}{P_0\left[\delta=0\vee\delta=\frac{\pi}{2}\vert\T\right]}-1\right]^{-1}
			\sim\left[1+\frac{\left[1-(p_0+p_{\frac{\pi}{2}})\right]\mathcal{I}_1}{\p f(\T|0)+\pp f(\T\big|\frac{\pi}{2})}\right]^{-1}
		\end{split}
	\end{equation*}
	as $\varepsilon\rightarrow0$.
	With the prior odds given by  
	\[
	R_0=\frac{p_0+p_{\frac{\pi}{2}}}{1-(p_0+p_{\frac{\pi}{2}})}
	\]
	and after algebraic manipulations, we obtain the approximation to the Bayes factor given by
	\begin{equation*}
		B_{01}
		\sim \frac{p_0f(\T|0)+p_{\frac{\pi}{2}}f(\T|\frac{\pi}{2})}{(p_0+p_{\frac{\pi}{2}})\mathcal{I}_1(\u{\theta})},
		~~\text{as $\varepsilon\rightarrow0$}.
	\end{equation*}
	\subsection{Test of vM axial symmetry}			\label{s23}
	
	We consider the Bayesian test of the null hypothesis that the sample follows a vM distribution against the 
	alternative that it comes from an arbitrary GvM distribution. This null hypothesis
	implies axial symmetry in the class of vM distributions, whereas the alternative 
	hypothesis includes both symmetric or asymmetric GvM distributions.  
	Precisely, we have
	$\text{H}_0:\kappa_2=0$,
	where $\mu_1,\mu_2$ and $\kappa_1$ are assumed known and equal to $\mu_1^0,\mu_2^0$ and $\kappa_1^0$ respectively. 
	The GvM with $\kappa_2 = 0$ reduces to the trivially symmetric vM distribution. 
	Formally, the GvM is defined for $\kappa_2>0$ only,
	so that the symmetry considered here is no longer within the GvM class but it is rather a vM axial 
	symmetry. This symmetry within the GvM class should be thought as  
	approximate, for vanishing values of $\kappa_2$. 
	
	Symmetry with the GvM formula can also be obtained with $\kappa_1=0$, in which case the GvM formula reduces to an axial von Mises distribution $\text{vM}_2$ that is trivially symmetric. This case is not analysed. In what follows we focus on the case of  vM axial symmetry. 
	
	Because $P[\kappa_2=0]=0$, we construct the perturbed probability $P_0$ such that
	$p_0=P_0[\kappa_2=0]>0$,
	where
	$
	p_0=P\left[\kappa_2\in\left[0,\varepsilon\right]\right]$, 
	for some $\varepsilon>0$ small.
	The prior d.f. of $\kappa_2$ under the probability $P$ is $G$,
and under the perturbed probability $P_0$ it is
	$
	p_0\Delta(\kappa_2')+(1-p_0)G(\kappa_2')$, $\forall\kappa_2' \ge 0$.
	Assume $0\notin(\kappa_2',\kappa_2'+d\kappa_2')$, then
	\[
	P_0\left[\kappa_2\in\left(\kappa_2',\kappa_2'+d\kappa_2'\right)\right]=(1-p_0)g(\kappa_2')d\kappa_2'=(1-p_0)P\left[\kappa_2\in\left(\kappa_2',\kappa_2'+d\kappa_2'\right)\right],
	\]
	where $g$ is the density of $G$.
	With algebraic manipulations similar to those of Section \ref{s21}, 
	one obtains the approximation to the Bayes factor $B_{01}$ of posterior over prior odds
	given by 
	\begin{equation}
		B_{01}\sim\frac{f(\u\theta|0)}{\int_{\mathcal{C}_{\varepsilon}^c}f(\u\theta|\kappa_2')g(\kappa_2')d\kappa_2'},~\text{as}~\varepsilon\rightarrow0,
		\label{MonteCarlo}
	\end{equation}
	where $\mathcal{C}_{\eps}=[0, \eps]$, 
	$\mathcal{C}_{\epsilon}^{\sf{c}}$ is its complement and where 
	the likelihood of $\kappa_2$ is
	\begin{align}    					                                \label{BFKappa2}
		f(\u{\theta}|\kappa_2) & =
		\left\{2\pi G_0\left(\delta^{0},\kappa_1^0,\kappa_2\right)\right\}^{-n}
		\exp\left\{\kappa_1^0 \sum_{i=1}^{n}\cos(\theta_i-\mu_1^0)
		+\kappa_2 \sum_{i=1}^{n} \cos 2 ( \theta_i - \mu_1^0 + \delta^0) \right\},
	\end{align}
	with $\delta^{0}=(\mu_1^{0}-\mu_2^{0})~\text{mod}~\pi$.
	\section{Numerical studies}							\label{s3} 
	This section provides some numerical studies for the tests introduced in Section \ref{s2}.
	The major part is Section \ref{s310}, which gives a simulation or Monte Carlo study  	
	of the performance of these tests. Section \ref{s40} provides an application to 
	real measurements of wind directions.  
	\subsection{Monte Carlo study}							\label{s310}
	This section presents a Monte Carlo study for
	the tests introduced in Section \ref{s2}:
	in Section \ref{s31} for the test of no shift between cosines,
	in Section \ref{s32} for the test axial symmetry and 
	in Section \ref{s33} for the test of vM axial symmetry.
	The results are summarized in Section \ref{s34}.
	We obtain Bayes factors for each one of these three tests for
	$r=10^{4}$ generations of samples of size $n=50$, that are generated from the GvM or the vM distributions. 
	The Monte Carlo approximation to the integral \eqref{e123}, and to the analogue integrals of the two other tests, is computed with $s=10^4$ generations.\\
	This simulation scheme is repeated three times and the results are compared in order to verify convergence. 
	Confidence intervals for the Bayes factors based on the aggregation of the three simulations (with $r$ replications each) are provided.
	
	The axial vM distribution ($\text{vM}_2$) 
	is used as a prior distribution for the parameter of
	shift between cosines $\delta$. This distribution can be obtained by taking $\kappa_1 = 0$ in
	the exponent of (\ref{e4}) and by multiplying the density by 2, yielding
	\begin{align}								\label{e321} 
		f(\theta \mid \mu,\kappa) & =
		\{\pi I_0(\kappa)\}^{-1} \exp \{\kappa \cos 2 (\theta -\mu)\}, \; 
		\; \forall \theta \in [0, \pi), 
	\end{align} 
	and for some $\mu \in [0, \pi)$ and $\kappa>0$.
	We denote an axial random variable with this distribution by
	$\text{vM}_2(\mu,\kappa)$.
	
	According to the remark at the end of Section \ref{s2new}, 
	we choose $\varepsilon=0.05$ for the length of the interval of $\text{H}_0$ and the prior densities $g$ as follows.
	For the test of no shift between cosines, we choose the ${\rm vM}_2(0,\tau)$ distribution 
	for $\delta$, which is symmetric and unimodal with mode at $\delta=0$. For the test of axial 
	symmetry, we choose the mixture of ${\rm vM}_2(0,\tau)$ and ${\rm vM}_2(\pi/2,\tau)$ for $\delta$.
	Finally, for the test of vM axial symmetry, we choose an 
	uniform distribution for $\kappa_2$ 
	that is highly concentrated at the boundary point $0$.
	
	\subsubsection{Test of no shift between cosines of GvM}	\label{s31} 
	
	The null hypothesis considered is $\text{H}_0$: $\delta=0$, with 
	fixed $\mu_1=\mu_1^0,\kappa_1=\kappa_1^0,\kappa_2=\kappa_2^0$, where 
	$\mu_1^0=\pi$, $\kappa_1^{0}=0.1$, $\kappa_2^0=5.5$. We consider three different cases, 
	called D1, D1' and D2.
	\\
	\textbf{Case D1}
	For $i=1,\ldots,s$,
	we generate $\delta^{(i)}$ from the prior of $\delta$, which is  
	$\text{vM}_2(\nu,\tau)$ with values of the hyperparameters $\nu=0$ and $\tau=250$.
	We obtain $p_0 = 0.570$ as prior probability of the null hypothesis under the perturbed 
	probability measure. 
	We take the first $r$ of these prior values (that are all the values, since $r=s$) and then we obtain $\mu_2^{(i)}=(\mu_1^0-\delta^{(i)})~\text{mod $\pi$}$
	and generate the elements of the vector of $n$ sample values 
	$\T^{(i)}$ independently from $\text{GvM}(\mu_1^{0}, \mu_2^{(i)},\kappa_1^{0},\kappa_2^{0})$, for $i=1,\ldots,r$.
	With these simulated data we compute the Bayes factor $B_{01}^{(i)}$ 
	with the approximation formula \eqref{BF_delta}.
	We repeat this experiment three times.
	The fact of generating values of $\delta$ from its prior distribution, 
	instead of taking $\delta =0$ fixed by null hypothesis, 
	is a way of inserting some prior uncertainty in the generated sample. 
	If the prior is close, in some sense, to the null hypothesis, 
	then we should obtain the Bayes factor larger than one, but
	smaller than the Bayes factor that would be obtained with the fixed value $\delta = 0$.   
	
\begin{figure} [h!]
		\begin{minipage}{.5\linewidth}
			\centering
			\subfloat[3 Boxplots of the 3 sets of $10^4$ simulated Bayes factors.]
			{\includegraphics[scale=0.25]{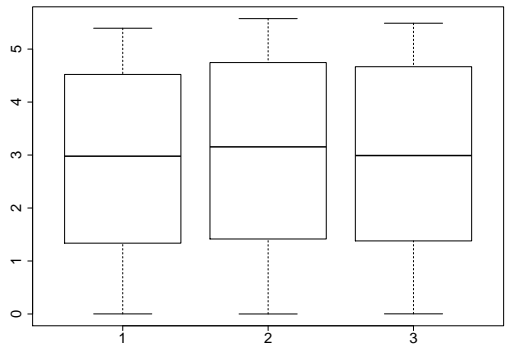}\label{delta_box}}
		\end{minipage}%
		\begin{minipage}{.5\linewidth}
			\centering
			\subfloat[Histogram of the sample of $3\cdot r$ Bayes factors and graph of its estimated density (red line).]
			{\includegraphics[scale=0.25]{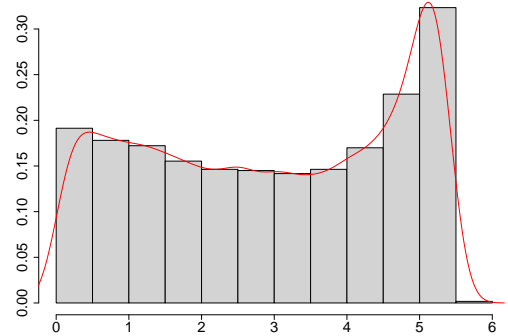}\label{normality_D1}}
		\end{minipage}\par\medskip
		\centering
		\subfloat[Prior density of $\delta$ with histogram of $10^4$ 
		generated values.]{\includegraphics[scale=0.25]{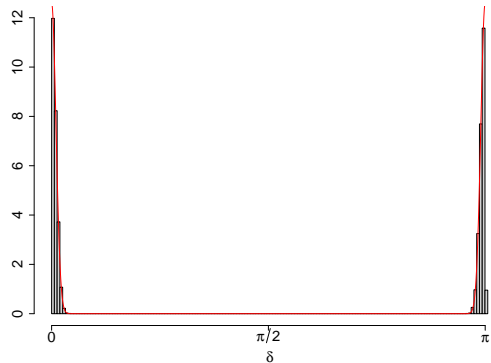}\label{delta_hist}}\quad
		\caption{\small Results of Case D1.}					
	\end{figure}

	We obtained three sequences of $10^4$ Bayes factors that can be summarized as follows.
	Figure \ref{delta_box} displays the three boxplots of the three simulated sequences of Bayes factors: 
	Denote by $\bar{B}_{01}^{(j)}$ the mean of the Bayes factors of the $j$-th sequence, for $j=1,2,3$,
	corresponding to left, central and right boxplot respectively.
	We obtained: 
	$$\bar{B}_{01}^{(1)} =2.887, \; \bar{B}_{01}^{(2)} = 3.028\; \text{ and } 
	\bar{B}_{01}^{(3)} =2.955. $$
	Figure \ref{normality_D1} shows the histogram of the three 
	generated sequences of $r$ 
	Bayes factors. The distribution is clearly not ``bell-shaped'' but it is however 
	light-tailed: the Central limit theorem applies to the mean of the simulated 
	Bayes factors. The asymptotic normal confidence interval for the mean 
	value of the Bayes factors at level $0.95$, and based on the 
	three generated sequences,
	is given by
	$$(2.937, 2.976).$$ 
	According to Table \ref{tab} this interval indicates positive evidence for the 
	null hypothesis: 
	the data have indeed increased the evidence of the null hypothesis that
	$\delta = 0$, however to a marginal extent only. 
	This situation can be explained by the fact that the prior density 
	$g$ is (highly) concentrated around $0$, circularly. This can be seen in the graph of the prior density (Figure \ref{delta_hist}), where
	the histogram of $10^4$ generated values of $\delta$ is shown together with the prior 
	density. Moreover, the variability originating from the fact the data are simulated under different values of $\delta$ leads to weaker values of the Bayes factor.
	\\
	\textbf{Case D1'} 
	In this other case we consider prior values of $\delta$ less concentrated around 0, 
	by choosing $\nu=0$ and $\tau=50$. The resulting prior probability of $\H_0$ is given by 
	$p_0 = 0.276$.
	For $i = 1, \ldots, r$, we generate the elements of the vector of $n$ sample values
	$\T^{(i)}$ independently from $\text{GvM}(\mu_1^{0}, \mu_2^{0},\kappa_1^{0},\kappa_2^{0})$,
	with $\delta = 0$, thus with $\mu_2^0 = (\mu_1^0 - \delta) \mod \pi = 0$.
	With these simulated data, we compute the Bayes factor $B_{01}^{(i)}$
	with the approximation formula \eqref{BF_delta}.
	
\begin{figure}[h!]
		\centering
		\subfloat[3 Boxplots of the 3 sets of $10^4$ simulated Bayes factors.]
		{\includegraphics[scale=0.3]{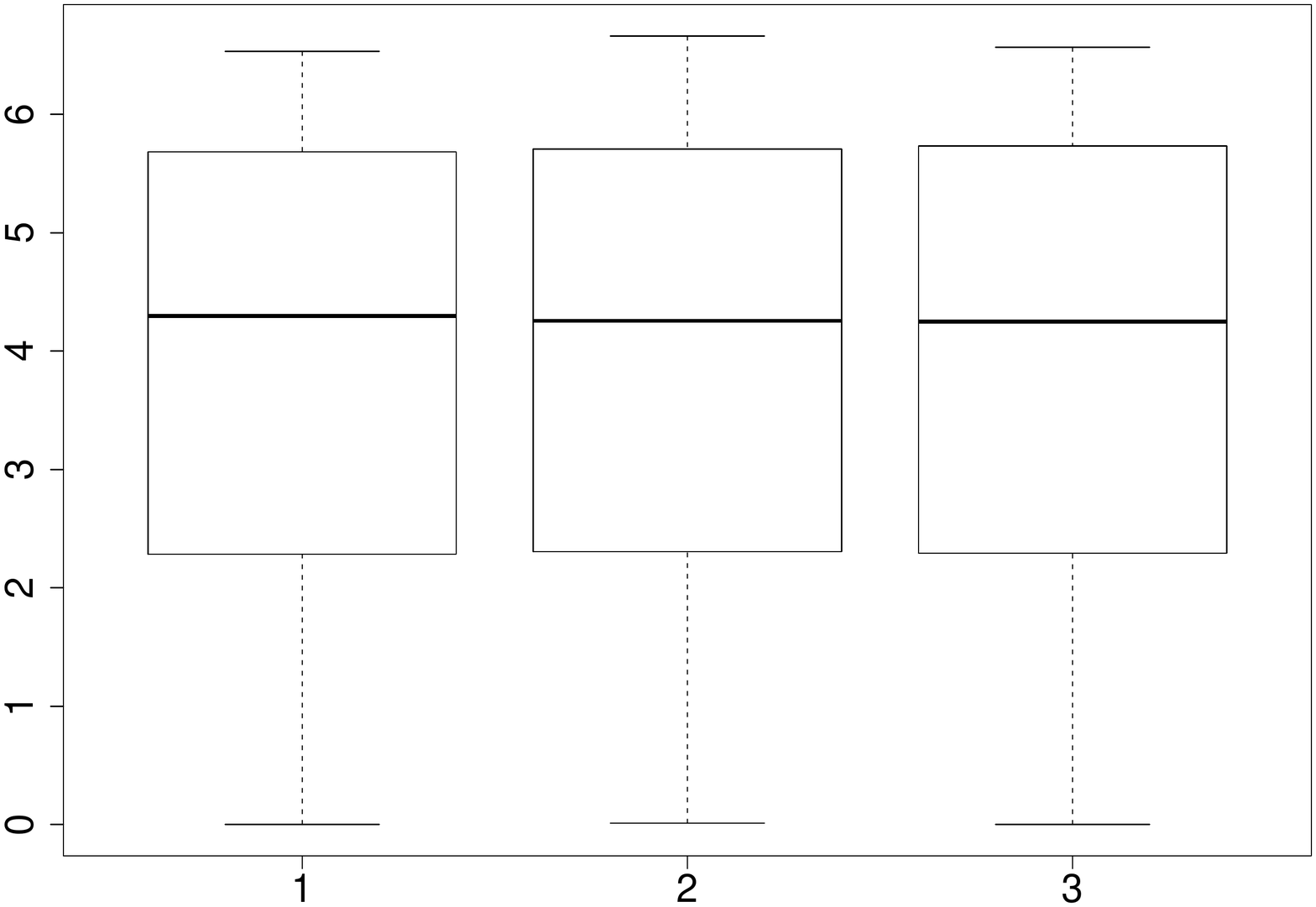}\label{delta_box'}}\\
		\subfloat[Prior density (dashed line) and posterior density (continuous line)
		of $\delta$. The posterior is based on one generated sample.]
		{\includegraphics[scale=0.3]{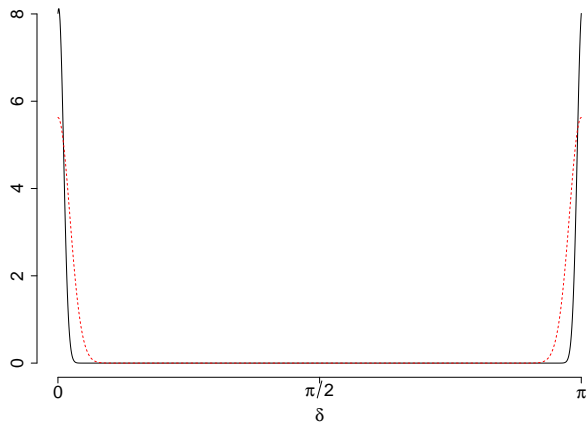}\label{D1'norm}}\\
		\caption{\small Results of Case D1'.}  
	\end{figure}

	We obtained three sequences of $r = 10^4$ Bayes factors with means: 
	$$\bar{B}_{01}^{(1)}=3.922, \; \bar{B}_{01}^{(2)}=3.924 \; \text{ and } 
	\bar{B}_{01}^{(3)}= 3.921.$$ 
	The boxplots of the three respective generated sequences are shown in Figure \ref{delta_box'}.
	The asymptotic normal confidence interval for the mean value of the Bayes factors, at level 0.95 and based on the 
	three generated sequences, is
	$$(3.901, 3.945).$$ 
	As expected, the generated Bayes factors 
	are larger than in case D1. 
	The samples
	generated with $\delta=0$ fixed have less uncertainty.
	We computed the posterior density of $\delta$ based on one generated sample. 
	In Figure \ref{D1'norm}
	we can see the graph of that posterior density, in continuous line, together with
	the graph of the prior density, in dashed line.
	The posterior is indeed more concentrated around 0, circularly.
	
	\noindent \textbf{Case D2} We now further decrease the concentration of the prior of $\delta$. The 
	values of the hyperparameters are $\nu=0$ and $\tau=20$.
	We computed the prior probability of the null hypothesis under perturbation $p_0 = 0.176$.
	We generated the samples 
	$\T^{(i)}$, for $i=1,\ldots,r$, with fixed 
	value $\mu_2^{0}=0$.
	
\begin{figure}[h!]
		\centering
		\subfloat[3 Boxplots of the 3 sets of $10^4$ simulated Bayes factors.]
		{\includegraphics[scale=0.3]{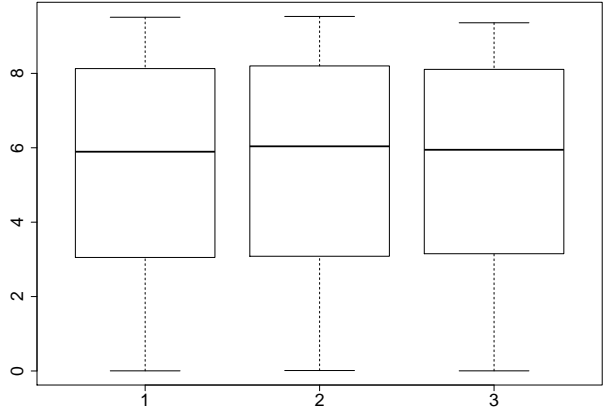}\label{boxplot_delta1}}\\
		\subfloat[Prior density (dashed line) and posterior density (continuous line)
		of $\delta$. The posterior is based on one generated sample.]
		{\includegraphics[scale=0.3]{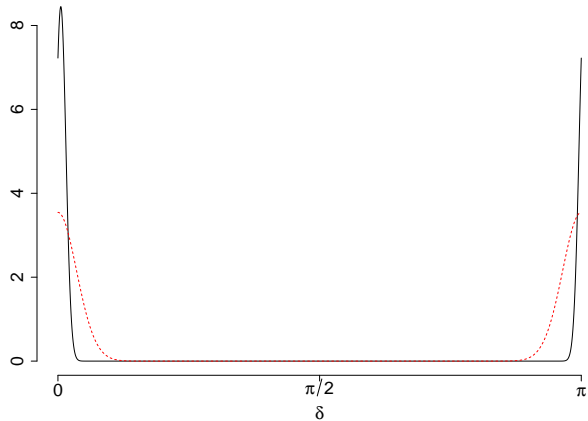}\label{delta_hist1}}\quad 
		\caption{\small Results of Case D2.}
	\end{figure}

	We obtained
	three sequences of $r = 10^4$ Bayes factors with means 
	$$\bar{B}_{01}^{(1)}=5.477,\bar{B}_{01}^{(2)}=5.539 \; \text{ and } \bar{B}_{01}^{(3)}=5.511.$$  
	The boxplots of the three respective generated sequences are shown in Figure \ref{boxplot_delta1}. 
	The asymptotic normal confidence interval for the mean value of the Bayes factors, at level 0.95 and based on the three 
	generated sequences, is
	$$(5.477, 5.541).$$ 
	The Bayes factors are larger than they are in Cases D1 and D1'.
	Here they show substantial evidence for the null hypothesis. 
	The prior distribution $\delta$ is less favourable to the null hypothesis and so the sample brings 
	more additional evidence for the null hypothesis.  
	Figure \ref{delta_hist1} shows the graph of the prior density, as dashed line, together with the 
	graph of a posterior density, as continuous line, for $\delta$.
	The graph of the posterior density is based on one generated sample.
	
	\subsubsection{Test of axial symmetry of GvM}		\label{s32}
	
\begin{figure} [h!]
		\centering
		\subfloat[3 Boxplots of the 3 sets of $10^4$ simulated Bayes factors.]
		{\includegraphics[scale=0.3]{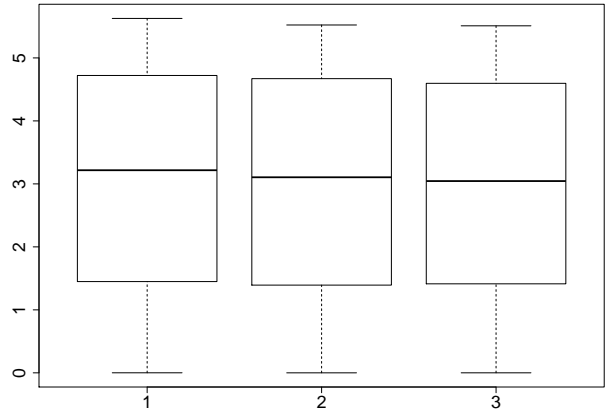}\label{symm_box}} \\
		\subfloat[Prior density of $\delta$ with histogram of $10^4$ generated values.]
		{\includegraphics[scale=0.3]{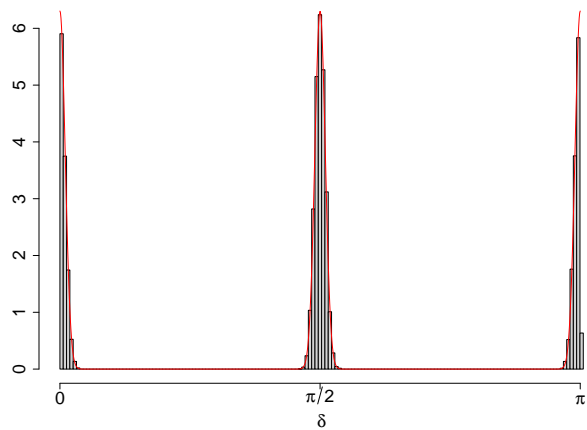}\label{symm_hist}}\quad
		\caption{\small Results of Case S1.}
	\end{figure}

	In this section we consider 
	the null hypothesis of axial symmetry, viz.
	$\text{H}_0$: $\delta=0$ or $\delta=\pi/2$,
	other parameters being fixed as follows,
	$\mu_1=\mu_1^0,\kappa_1=\kappa_1^0$ and $\kappa_2=\kappa_2^0$. We choose
	as before $\mu_1^0=\pi$, $\kappa_1^0=0.1$ and $\kappa_2^0=5.5$. 
	We generate $\delta$ from the prior given by the
	mixture of $\text{vM}_2$ distributions
	$\xi \, \text{vM}_2(\nu_1,\tau)+(1-\xi) \, \text{vM}_2(\nu_2,\tau)$, 
	with $\nu_1=0,\nu_2=\pi/2$ and $\xi = 0.5$.
	We consider three different cases, called Cases S1, S2 and S3.
	\\
	\textbf{Case S1} We generated $\delta$ from the prior mixture 
	with concentration parameter $\tau=250$. 
	This prior distribution is close to the null distribution and 
	Figure \ref{symm_hist} displays its density, 
	together with the histogram of $10^4$ generations from it.
	We computed the prior probabilities of the null hypothesis under the 
	perturbed probability measure $p_0 = p_{\pi/2} = 0.285$.
	We follow the principle of Case D1, where prior uncertainty is
	transmitted to the sample by considering generated values $\delta^{(i)}$, for $i=1,\ldots,s$, 
	from a prior of $\delta$ close to the null hypothesis, instead of considering
	the fixed values of the null hypothesis, namely $\delta = 0$ or $\pi/2$. We take the first $r$ of these prior values and we use $\mu_2^{(i)}=(\mu_1^{0}-\delta^{(i)})~\text{mod }\pi$ for generating 
	$\T^{(i)}$, for $i=1,\ldots,r$. Repeating this three times, we obtained 
	the three means of the three sequences of $r=10^4$ Bayes factors 
	$$\bar{B}_{01}^{(1)}=3.044,\bar{B}_{01}^{(2)}=2.986 \text{ and } \bar{B}_{01}^{(3)}=2.950.$$
	In Figure \ref{symm_box} we can find the boxplots 
	of the three respective generated sequences.
	The asymptotic normal confidence interval for the mean value of the  Bayes factors,  at level 0.95 and based on the three generated sequences, is
	$$(2.974,3.013).$$
	The conclusion is that the sample provides positive evidence of axial symmetry, 
	even though to some smaller extent only. 
	The same was found in Case D1.
	\begin{figure}[h!]
		\centering
		\subfloat[3 Boxplots of the 3 sets of $10^4$ simulated Bayes factors.]
		{\includegraphics[scale=0.3]{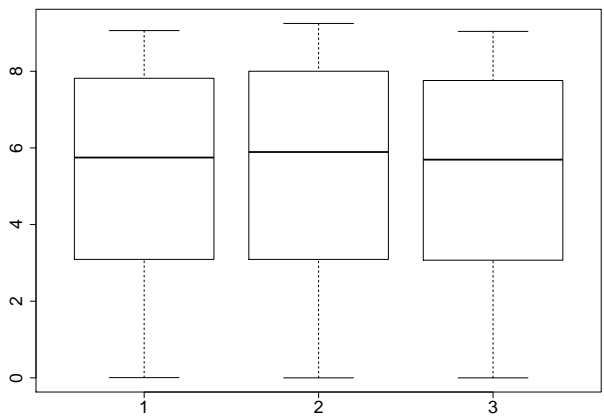}\label{symm_box1}}\\
		\subfloat[Prior density (dashed line) and posterior density (continuous line) of 
		$\delta$. The posterior is based on one generated sample.]{\includegraphics[scale=0.3]
			{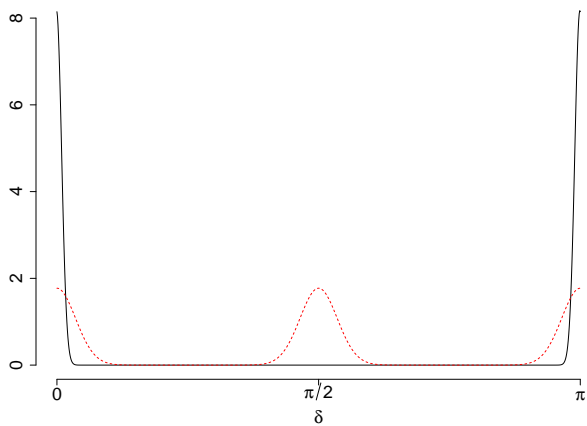}\label{symm_hist1}}\quad 
		\caption{\small Results of Case S2.}
	\end{figure}
	\begin{figure} [h!]
		\centering
		\subfloat[3 Boxplots of the 3 sets of $10^4$ simulated Bayes factors.]
		{\includegraphics[scale=0.3]{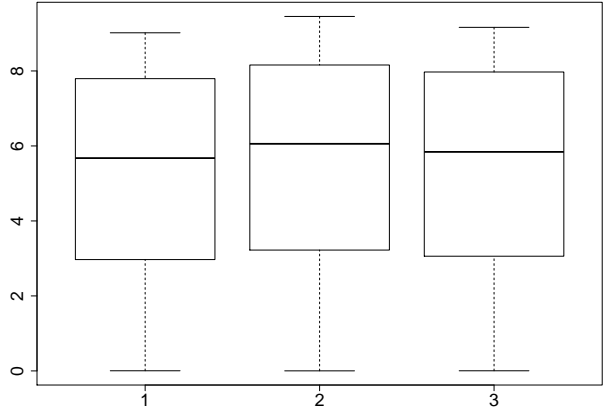}\label{symm_box2}} \\
		\subfloat[Prior density (dashed line) and posterior density (continuous line) of 
		$\delta$. The posterior is based on one generated sample.]
		{\includegraphics[scale=0.3]{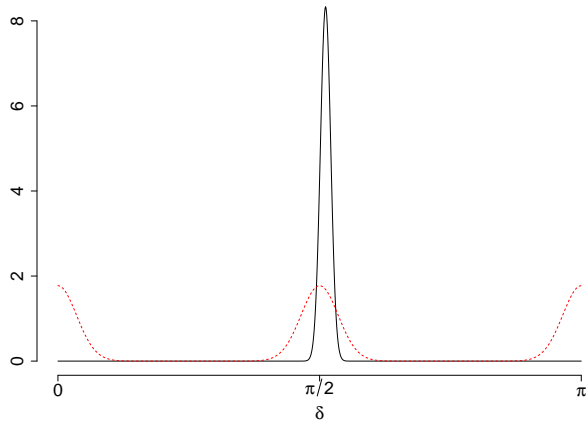}\label{S3B}}
		\caption{\small Results of Case S3.}
	\end{figure}\\
	\textbf{Case S2}  
	We generated prior values of $\delta$ from the same mixture, however with smaller concentration hyperparameter $\tau=20$. We found $p_0 = p_{\pi/2} = 0.088$.
	We generated the elements of the sample vector $\T^{(i)}$ 
	with fixed value $\mu_2^{0}=0$, thus from $\text{GvM}(\mu_1^{0}, \mu_2^{0},\kappa_1^{0},\kappa_2^{0})$, 
	with $\mu_1^{0}=\pi,\mu_2^{0}=0,\kappa_1^{0}=0.1,\kappa_2^{0}=5.5$, for $i=1,\ldots,r$. 
	We repeated this experiment three times and obtained 
	three sequences of Bayes factors, with respective mean values
	$$\bar{B}_{01}^{(1)}=5.322, \bar{B}_{01}^{(2)}=5.439 \; 
	\text{ and }\bar{B}_{01}^{(3)}=5.282.$$ The boxplots of the three sequences of 
	Bayes factors can be found in Figure \ref{symm_box1}.
	After aggregating the three sequences, we obtained the asymptotic normal confidence 
	interval at level $0.95$ for the mean value of the
	Bayes factors given by  
	$$(5.317, 5.378).$$
	The Bayes factor is thus larger than it was in Case S1, 
	so that the sample has brought substantial evidence of axial symmetry.
	Figure (\ref{symm_hist1}) shows the prior density of $\delta$ (dashed line) 
	and a posterior density of $\delta$ (continuous line) that is 
	based on one of the previously generated samples.
	The posterior is highly concentrated around $0$ 
	and provides a stronger belief about symmetry than the prior.
	\\ 
	\textbf{Case S3} We retain the prior of $\delta$ of Case S2 but we generate 
	samples $\T^{(i)}$, for $i = 1, \ldots, r$, with
	$\mu_1^{0}=\pi, \mu_2^{0}=\pi/2,\kappa_1^{0}=0.1$, and $\kappa_2^{0}=5.5$,
	thus from another symmetric GvM distribution.
	The computed values $p_0 = p_{\pi/2} =  0.088$ are the same of Case S2. 
	We generated three sequences of $r=10^4$ Bayes factors.
	The three respective boxplots of the three sequences can be found in \ref{symm_box2}.
	The three respective means of these three sequences are  
	$$\bar{B}_{01}^{(1)} = 5.267, \bar{B}_{01}^{(2)}=5.553 \text{ and } \bar{B}_{01}^{(3)}=5.395.$$
	By aggregating the three sequences, we obtained the asymptotic normal confidence
	interval at level $0.95$ for the mean of the 
	Bayes factors given by
	$$(5.374, 5.436).$$
	We find substantial evidence of axial symmetry.
	Figure \ref{S3B} displays the prior density of $\delta$ (dashed line)
	and a posterior density of $\delta$ (continuous line) that is
	based on one of the previously generated samples.
	The posterior is highly concentrated around $\pi/2$
	and possesses less uncertainty about symmetry than the prior.
	
	\subsubsection{Test of vM axial symmetry of GvM}   		\label{s33}
	
\begin{figure}[h!]
		\centering
		\includegraphics[scale=0.3]{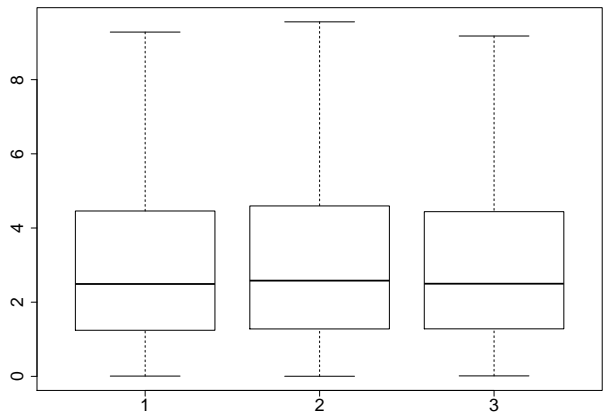}
		\caption{3 boxplots of the 3 sets of $10^4$ simulated Bayes factors.}
		\label{Boxplot_kappa2}
	\end{figure}

	Now we have 
	$\text{H}_0:\kappa_2=0$, with fixed $\mu_1^{0}=\pi$, $\mu_2^{0}=\pi/2$ and $\kappa_1^{0}=0.1$. The prior distribution of $\kappa_2$ is uniform over $[0,1/2]$ and the sample $\u\theta=(\theta_1,\ldots,\theta_n)$ is generated from the $\text{vM}(\mu_1^{0}, \kappa_1^{0})$ distribution. The prior probability of $\text{H}_0$ under the perturbation is $p_0=0.1$. We generated three sequences of $r=10^4$ Bayes factors: their boxplots are shown in Figure \ref{Boxplot_kappa2}. In these boxplots we removed a very small number of large values, in order to improve the readability. The three means of the three generated sequences are
	\[
	B_{01}^{(1)}=3.284,~B_{01}^{(2)}=3.380 ~ \text{ and } B_{01}^{(3)}=3.241,
	\]
	where the very large values that were eliminated from the boxplots have been considered in the calculations of these means. 
	After aggregating these three sequences, we obtained the following asymptotic normal confidence interval for the mean value of the Bayes factors at level $0.95$, 
	\[
	(3.268, 3.335 ).
	\]
	There is a positive evidence of symmetry although rather limited. 
	The amount of evidence is similar to the cases D1 and S1: in all these studies, the prior is much concentrated around the null hypothesis (here $\kappa_2=0$), 
	so that the data have increased the evidence of the null hypothesis only to some limited extend. 
	
	\subsubsection{Summary}									\label{s34} 
	
	Table \ref{tab2} summarizes the simulation
	results that we obtained for the three tests and for the various cases.
	\begin{table}[h!]
		\centering
		{\small
			\begin{tabular}{|c|c|c|c|}
				\hline
				$\H_0$ & case & confidence interval for $B_{01}$ & evidence for $\H_{0}$\\
				\hline
				\multirow{3}{*}{no shift between cosines} & D1 & (2.937, 2.976)& positive\\
				&D1' & (3.901, 3.945) & positive\\
				& D2 & (5.477, 5.541) & substantial\\
				\hline
				\multirow{3}{*}{axial symmetry} & S1& (2.974, 3.013) & positive \\
				& S2 & (5.317, 5.378) & substantial\\
				& S3 & (5.374, 5.436) & substantial\\
				\hline
				vM axial symmetry & -- & (3.268, 3.335) & positive\\
				\hline
			\end{tabular}
			\caption{\small Summary of the simulation study.}
			\label{tab2}
		}
	\end{table}
	
	\subsection{Application to real data}				\label{s40} 
	
	The proposed Bayesian tests have been so far applied to simulated data. 
	This section provides the application of the test of no shift 
	between cosines of Section \ref{s21} and of 
	axial symmetry of Section \ref{s22} to real data obtained from 
	the study ``ArticRIMS'' (A Regional, Integrated Hydrological Monitoring 
	System for the Pan Arctic Land Mass) available at \url{http://rims.unh.edu}. 
	The Arctic climate, its vulnerability, its
	relation with the terrestrial biosphere and with the recent 
	global climate change are the subjects under investigation.
	For this purpose, various meteorological variables such as temperature, precipitation, 
	humidity, radiation, vapour pressure, speed and directions of winds 
	are measured at four different sites.
	
\begin{figure}[h!]
		\centering
		\includegraphics[scale=0.3]{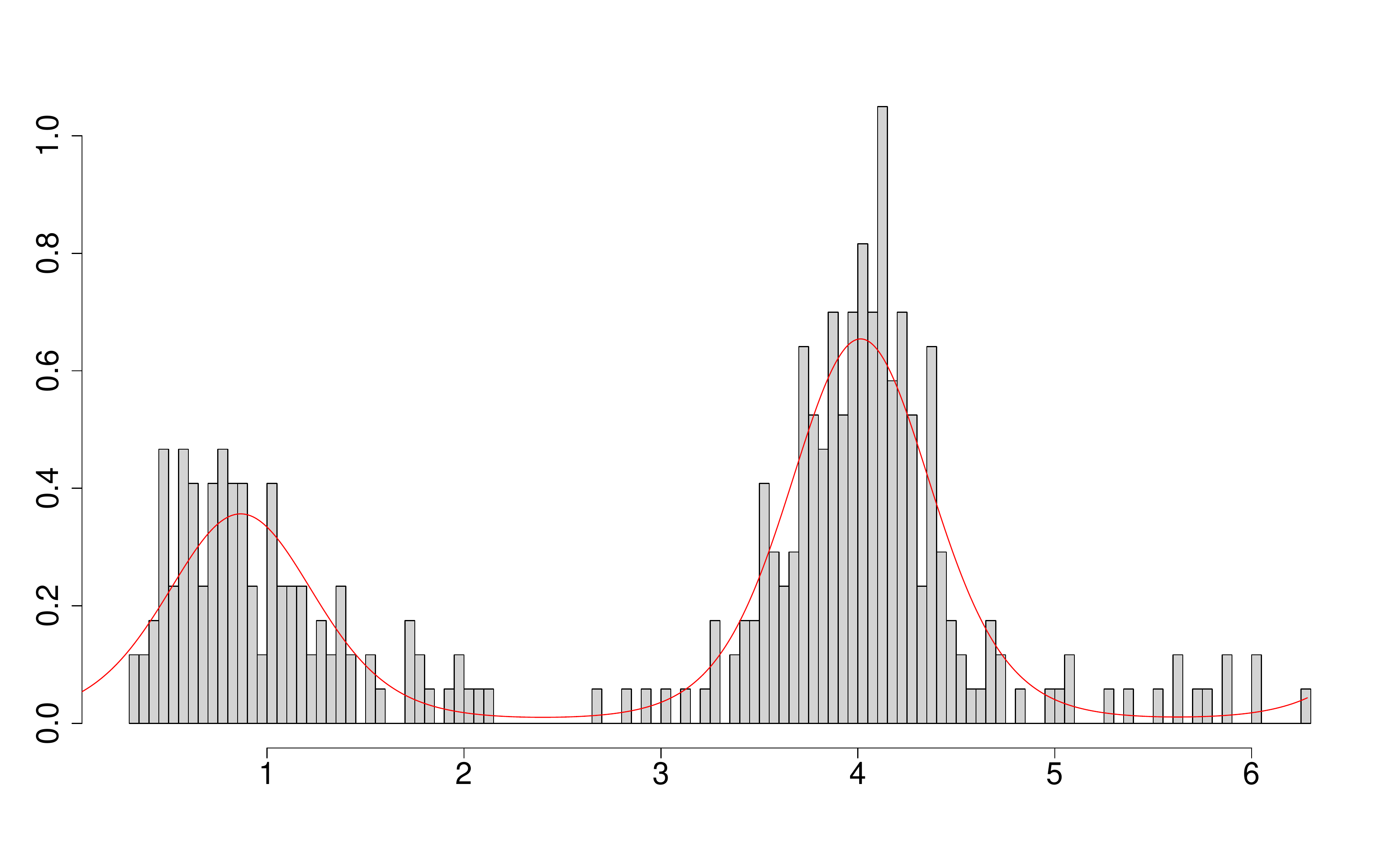}
		\caption{Histogram of wind directions data together with
			$\text{GvM}(\hat{\mu}_1,\hat{\mu}_2,\hat{\kappa}_1,\hat{\kappa}_2)$ density.}
		\label{europe05}
	\end{figure}

	We consider wind directions measured at the site 
	``Europe basin'' and from January to December 2005.
	After removal of few influential 
	measurements, 
	the following maximum likelihood estimators are obtained: 
	$\hat{\mu_1}=4.095,~~\hat{\mu}_2=0.869,~~\hat{\kappa}_1= 0.304,~~\hat{\kappa}_2=1.910$ and
	thus $\hat{\delta}=(\hat{\mu_1}-\hat{\mu}_2)\text{ mod }\pi=0.084$.
	The histogram of the sample together with the GvM density with theses values of the parameters 
	are given in Figure \ref{europe05}. 
	
	For the test of no shift between cosines, 
	the Monte Carlo integral (\ref{e123}) is computed with 
	$s=10^6$ values of $\delta$ generated from the prior $\text{vM}_2(\nu, \tau)$, 
	with $\nu=0$ and $\tau=300$. 
	We consider $\varepsilon=0.18$:
	as mentioned in Section \ref{s2new}, 
	a substantial value is desirable in the practice.
	We obtain the Bayes factor $B_{01} = 2.550$; cf. Table \ref{table_europe}.
	
	For the test of symmetry,  
	the prior of $\delta$ is the mixture of two vM of order two, 
	i.e. $\xi \, \text{vM}_2(\nu_1,\tau)+(1-\xi) \, \text{vM}_2(\nu_2,\tau)$, 
	with $\nu_1=0,~\nu_2=\pi/2,~\tau=300$ and $\xi=0.5$. 
	Monte Carlo integration 
	is done with $s=10^6$ generations from this prior.
	We consider  $\varepsilon=0.18$ and 
	obtain the Bayes factor $B_{01} = 2.252$; cf. Table \ref{table_europe}.
	
	The values of the two Bayes factors of Table \ref{table_europe} show positive
	evidence for the respective null hypotheses. 
	\begin{table}[h!]
		\centering
		{\small 
			\begin{tabular}{|c|c|c|}
				\hline
				$\text{H}_0$ & $B_{01}$& evidence for $\text{H}_{0}$\\
				\hline
				no shift between cosines & 2.550 & positive\\
				axial symmetry & 2.252 & positive\\
				\hline
			\end{tabular}
			\caption{Bayes factors for wind directions data.}
			\label{table_europe}
		}
	\end{table}
	
	\section{Conclusion}								\label{s4} 
	
	This article introduces three Bayesian tests relating to the symmetry of the GvM model.
	The first test is about the significance of the shift parameter 
	between the cosines of frequency one and two ($\H_0: \delta = 0$).
	The second test is about axial symmetry
	($\H_0: \delta = 0$ or $\delta = \pi/2$). The third test is about vM symmetry ($\text{H}_0:\kappa_2=0$).
	These tests are obtained by the technique of probability perturbation.
	Simulation studies show the effectiveness of these three tests, in the sense
	that when the sample is coherent with the null hypothesis, then 
	the Bayes factors are typically large. 
	Applications to real data are also shown.
	
	Due to computational limitations, we consider 
	null hypotheses of symmetry that concern one parameter only.
	The null hypotheses considered are about one or two distinct values 
	of the parameter of interest, with all remaining parameters fixed.
	Composite null hypotheses that allow for unknown nuisance parameters, would 
	require one additional dimension of Monte Carlo integration for each unknown parameter,
	in the computation of the marginal
	distribution.
	The computational burden would rise substantially and the Monte Carlo study,
	with two levels of nested generations, would become very difficult.
	But the essentially simple null hypotheses considered
	are relevant in the practice. It can happen that nuisance parameters
	have been accurately estimated and the question of interest is really about the
	the parameter $\delta$ and axial
	symmetry. In the example of Section \ref{s40}, we want to know if 
	wind direction is axially symmetric within the
	GvM model. The values of the concentrations and of the 
	axial direction are of secondary importance.
	
	One could derive other Bayesian tests for the GvM model: a 
	Bayesian test of bimodality is under investigation.
	We can also note that \cite{AAAI1715020} introduced an useful multivariate 
	GvM distribution for which similar Bayesian tests could be investigated.
	
	The computations of this article are done with the language R, see 
	\cite{r-project}, over a computing cluster with several cores. 
	The programs are available at the software section of 
	{\tt http://www.stat.unibe.ch}. 
	
	\appendix 
	
	\section{Proof of Proposition \ref{p1}}						\label{a1}   
	
	The definition of axial symmetry given at the beginning of Section 
	\ref{s22} tells that
	the GvM distribution is symmetric around $\alpha/2$ (or $\alpha/2 + \pi$),
	for some $\alpha \in [0 , 2 \pi)$, iff
	\begin{equation*}
		f(\theta|\mu_1,\mu_2,\kappa_1,\kappa_2)=f(\alpha-\theta|\mu_1,\mu_2,\kappa_1,\kappa_2),
		~~\forall\theta\in[0,2\pi).
		\label{symmetry}
	\end{equation*}
	This means
	\begin{equation}
		\begin{split}
			\kappa_1\cos(\theta-\mu_1)+\kappa_2\cos2(\theta-\mu_2)&=\kappa_1\cos[(\alpha-\theta)-\mu_1]+\kappa_2\cos2[(\alpha-\theta)-\mu_2]\\
			&=\kappa_1\cos[\theta-(\alpha-\mu_1)]+\kappa_2\cos2[\theta-(\alpha-\mu_2)],
		\end{split}
		\label{ax}
	\end{equation}
	$\forall \theta\in[0,2\pi)$.
	By using the cosine addition formula, \eqref{ax} can be re-expressed as
	\begin{equation*}
		\begin{split}
			\kappa_1\cos\theta\cos\mu_1+\kappa_1\sin\theta\sin\mu_1+\kappa_2\cos2\theta\cos2\mu_2+\kappa_2\sin2\theta\sin2\mu_2 = \qquad  \qquad \qquad \qquad \qquad  \\
			\kappa_1\cos\theta\cos(\alpha-\mu_1)+\kappa_1\sin\theta\sin(\alpha-\mu_1)+\kappa_2\cos2\theta\cos2(\alpha-\mu_2)+\kappa_2\sin2\theta\sin2(\alpha-\mu_2),
		\end{split}
	\end{equation*}
	$\forall\theta\in[0,2\pi)$.
	This is equivalent to the equation
	\begin{equation*}
		\begin{split}
			&\kappa_1[\cos\mu_1-\cos(\alpha-\mu_1)]\cos\theta+\kappa_2[\cos2\mu_2-\cos2(\alpha-\mu_2)]\cos2\theta\\
			&\qquad+\kappa_1[\sin\mu_1-\sin(\alpha-\mu_1)]\sin\theta
			+\kappa_2[\sin2\mu_2-\sin2(\alpha-\mu_2)]\sin2\theta=0,
		\end{split}
	\end{equation*}
	$\forall\theta\in[0,2\pi)$.
	It is convenient to re-express this last equation
	in terms of a trigonometric polynomial of degree $N=2$, precisely as
	\begin{equation}
		p(\theta)=\sum_{j=1}^{N}(a_j \cos j \theta + b_j \sin j \theta)=0,~\forall\theta\in [0,2\pi),
		\label{p}
	\end{equation}
	whose coefficients are given by
	\begin{align*}
		a_j & =\kappa_j[\cos j\mu_j-\cos j(\alpha-\mu_j)] \; \text{ and } \;
		b_j = \kappa_j[\sin j\mu_j-\sin j(\alpha-\mu_j)],
		\text{ for } j=1,2.
	\end{align*}
	A trigonometric polynomial of degree $N$ has maximum $2N$ roots in $[0,2\pi)$,
	unless it is the null polynomial; see e.g. p. 150 of \cite{powell1981approximation}.
	With this, \eqref{p} implies that $p(\theta)$ is the null polynomial,
	which means that $a_j = b_j = 0$, for $j=1,2$. These four equalities give
	the system of equations
	\begin{equation*}
		\begin{cases}
			\cos\mu_1=\cos(\alpha-\mu_1), \\
			\sin\mu_1=\sin(\alpha-\mu_1), \\
			\cos2\mu_2=\cos2(\alpha-\mu_2), \\
			\sin2\mu_2=\sin2(\alpha-\mu_2),
		\end{cases}
	\end{equation*}
	which, in terms of $\delta = (\mu_1 - \mu_2) \mod \pi$,
	simplifies to
	\begin{equation*}
		\begin{cases}
			\alpha=2\mu_1+2k\pi, \\
			\alpha=2\big[(\mu_1-\delta)~\text{mod $\pi$}\big]+k_1\pi,
		\end{cases} \; \text{ for some } k,k_1\in\Z.
	\end{equation*}
	One can eliminate the congruence symbol mod and obtain
	\begin{equation}                                                                \label{sim}
		\begin{cases}
			\alpha = 2\mu_1+2k\pi, \\
			\alpha = 2\mu_1-2\delta + (k_1+2k_2) \pi ,
		\end{cases}
		\; \text{ for some } k,k_1,k_2\in\Z.
	\end{equation}
	This system of simultaneous equation admits solutions iff $2\delta$ is a multiple of $\pi$,
	i.e.  
	$2\delta=0~\text{mod $\pi$}$.
	Since $\delta\in[0,\pi)$, we have found the desired symmetry characterization.
	
	{\small 
		\bibliography{manuscript}
		\bibliographystyle{apalike}
		\setcitestyle{authoryear,open={(},close={)}}
	}
	
\end{document}